\documentclass[12pt]{article}
\usepackage{epsfig}
\usepackage{amssymb}
\title{Build Boy's Surface}
\author{Richard Evan Schwartz \thanks{Supported by
    N.S.F. Grant DMS-2505281}}

\def\R{\mbox{\boldmath{$R$}}}%

\begin{document}
\maketitle

\section{Introduction}

Boy's surface [{\bf B\/}] is a famous immersed copy of
the projective plane in $\R^3$.
I have tried at various times to
understand Boy's surface, but somehow the
effort has never met with joy.
From time to time I have been able to concentrate
enough so that its crazy twists and turns 
 fit together in my mind (or seemed to), but then rather quickly
 the image went away without a trace.
 Moira Chas recently gave me one of her
 beautiful wire-knit models [{\bf C\/}] of Boy's surface,
 and this inspired me to try again.
 Now I have it straight. 
 For a nice and succinct
 alternate description of Boy's surface, see [{\bf K\/}].

These notes give a fairly conceptual
description of Boy's surface that does not draw too much
on three dimensional visualization.
  As I discuss in \S \ref{airplane}, 
 the approach I take owes a big intellectual debt to
 a beautiful youtube video [{\bf S\/}],
 produced by the Serbian
  Academy of Sciences, showing an airplane flying
  around and constructing Boy's surface.
  You could view these notes as an elaboration
  of the video.

I will
explain how to build Boy's surface out of
simple pieces.  All the pieces except the last one are easy to make
out of cardstock and tape, and the last piece is determined
automatically from the others by a coning procedure.
For the coning procedure to work most gracefully it is
useful to think of Boy's surface as living in the
$3$-sphere, which I model as $\R^3 \cup \infty$.

I will include a kit that lets you print
out and assemble all but the last piece.
The kit works.
Katherine Williams Booth made a model from it.
Here is her photograph of her model.
    
\begin{center}
  \resizebox{!}{5.2in}{\includegraphics{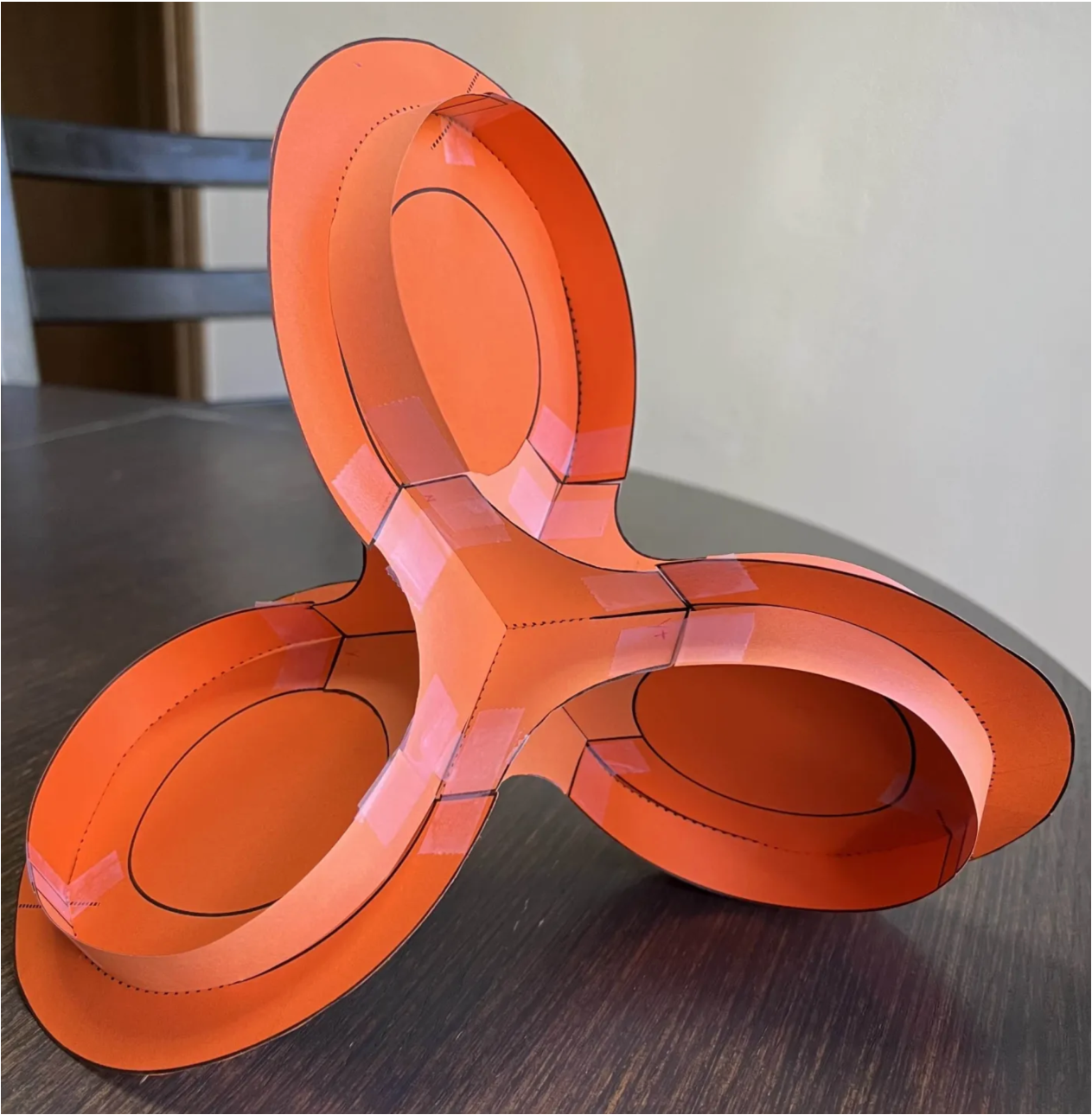}}
  \newline
  Figure 0:  An assembled kit!
\end{center}

Here is an outline of these notes.
\newline
$\bullet$
  In \S 2, I discuss the {\it octahedral graph\/} in some
  detail.   This discussion gives a combinatorial explanation of why
  the
  boundary of the immersed surface in Figure 0 is
  a single loop.  Of course, you can also just build the
  thing and trace the boundary around with your finger.
  \newline
  $\bullet$
  In \S 3, I define the {\it octacross\/}.
  This piece is homeomorphic to the union of the
  $3$ coordinate planes.  It is the central piece in Figure 0. \newline
  $\bullet$
  In \S 4, I define the {\it crossbridge\/}.
  This piece is homeomorphic to the union of $2$ coordinate
  planes.  The $3$ crossbridges are the big hoops in Figure 0.
\newline
  $\bullet$ In \S 5, I assemble most of Boy's
  surface:  You attach $3$ crossbridges to the octacross
  and then glue in $3$ round disks.  The disks are
  visible in Figure 0.  The union $M$ of these
  $7$ pieces turns out to be an immersed Moebius band.
  This is what Figure 0 shows.
  \newline
  $\bullet$
  In \S 6, I attach the final piece $P$, which I call
  {\it the pizza\/}.   It turns out that the boundary
  $\partial M$ is {\it cone-friendly\/}:  Each ray through
  the origin intersects $\partial M$ at most once.
  We set $P=\{tp|\ p \in \partial M,\ t \geq 1\}$.
  Here $P$ is a topological disk we get by coning
  $\partial M$ to $\infty$.
  The reason for the name is that $\partial M$ 
  has the combinatorial structure of an $18$-gon,
  and correspondingly $P$ is made from
  $18$ {\it pizza slices\/}, each one a topological triangle.
  Boy's surface is $M \cup P$.
  \newline
  $\bullet$
  In \S 7, I  compare and contrast these notes with the
  youtube video [{\bf S\/}].  I also make a connection to paper
  Moebius bands [{\bf Sch\/}] and raise an optimization question
related to the construction here. \newline
  $\bullet$
  In \S 8, I include the kit which lets you make $M$, as in Figure 0.
  \newline
  $\bullet$ In \S 9, I present an alternate model
  using rectilinear pieces.

I   thank Moira Chas, Peter Doyle, and Dan Margalit for
   helpful conversations about these notes. Finally,
   I thank Katherine Williams Booth for making the model in Figure 0.
  
  \section{The Octahedral Graph}
  \label{octa}

  Let $G$ be the {\it octahedral graph\/}.
  This graph is the union of vertices and edges of the octahedron.
  A {\it perfect matching\/} for $G$ is a union
  of $3$ edges of $G$ such that every vertex of $G$ is contained in
  exactly one edge in the union.   Figure 1 shows a perfect matching
  for $G$, drawn in a symmetric way.
  
\begin{center}
  \resizebox{!}{2in}{\includegraphics{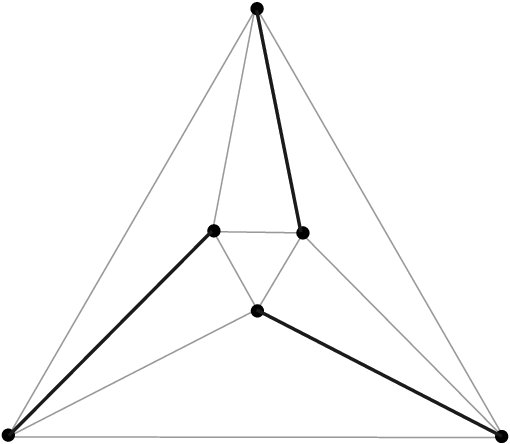}}
  \newline
  Figure 1:  A perfect matching for the octahedral graph
  \end{center}

  We call the $3$ edges in the matching {\it bridges\/},
  and the remaining $9$ edges {\it roads\/}.
  The reason for the name is that each bridge $e$
  connects the roads that are incident to vertices
  on either side of $e$.
  We can thicken $e$ slightly and then replace $e$ with three
  parallel copies. (I like to think of these parallel copies as
  lanes on the bridge.)
  This defines $3$ paths of length $3$,
  each of which has the form road-bridge-road.
  Figure 2 shows this.
  
\begin{center}
  \resizebox{!}{1in}{\includegraphics{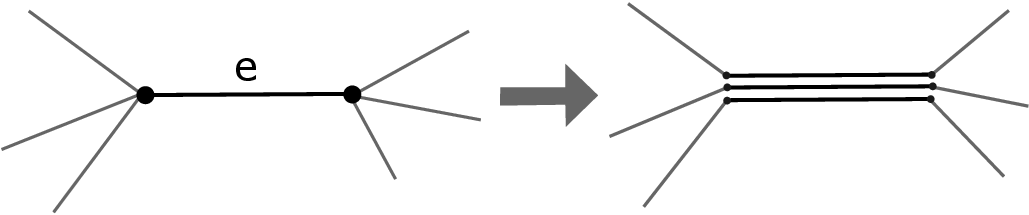}}
  \newline
  Figure 2:  A bridge between incident roads
\end{center}

Figure 3 shows what happens when we do this construction
simultaneously for all three bridges.  The result is a
length $18$ circuit that alternates between roads and
bridges.  I have added numbers and directions to the roads to help
you trace out the circuit.

\begin{center}
  \resizebox{!}{3.7in}{\includegraphics{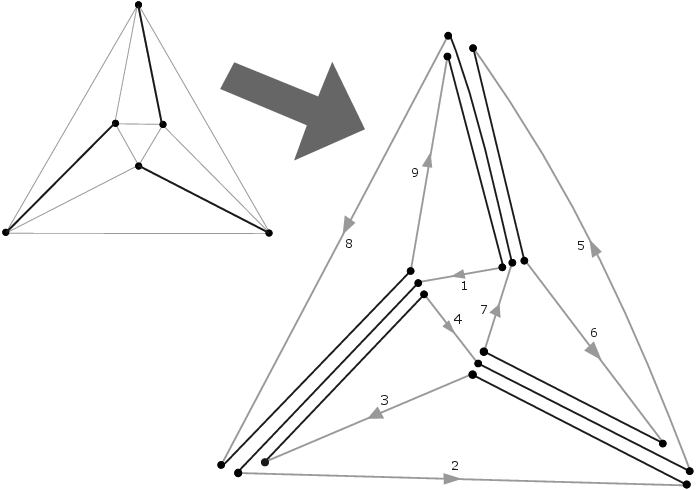}}
  \newline
  Figure 3:  A circuit of length $18$
    \end{center}

    Figure 4 is the same as Figure 3, but with different labels.
    The labels correspond to the usual choice of octahedron having
    vertices $(\pm 1,0,0)$ and $(0,\pm 1,0)$ and $(0,0,\pm 1)$. 
Geometrically, you are staring at the octahedron along the axis
$X=Y=Z$.  The vertex $(\pm 1,0,0)$ is labeled by $(\pm, 0,0)$. etc.
We label the roads by the signs of the coordinates of points on them.

\begin{center}
  \resizebox{!}{3.7in}{\includegraphics{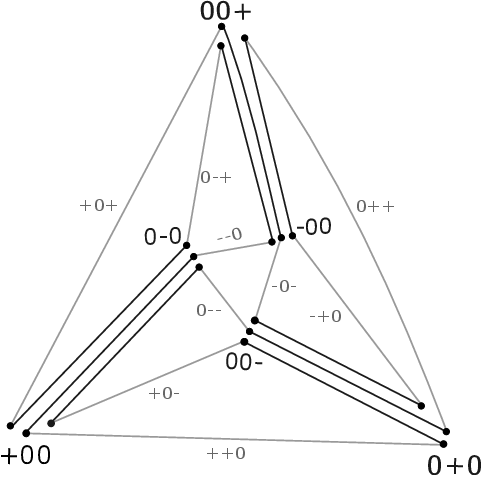}}
  \newline
  Figure 4:  Same path as in Figure 3, with different labels
  \end{center}

All our constructions below are invariant under the order
$3$ rotation about the axis $X=Y=Z$, namely
the map $(X,Y,Z) \to (Y,Z,X)$.
 Likewise, all our constructions
are invariant under the order $2$ rotations that respectively
fix the lines $X+Y=0$ and $Y+Z=0$ and $Z+X=0$.
The first of these rotations is given by the map $(X,Y,Z) \to
(-Y,-X,-Z)$.
The labelings in Figure 4 reflect these symmetries.

Referring to Figure 0, the $18$ vertices in Figure 4 correspond
to the $18$ points of non-smoothness on the surface boundary -- i.e.,
the $18$ points where the octacross attaches to the crossbridges.
Most of the tape appears around these points.

\section{The Octacross}

We first mention some notation which we will use repeatedly.
\begin{itemize}
\item Given a set $A_{XY}$, we get $A_{YZ}$ and $A_{ZX}$ by
    cycling the coordinates.
    \item $D_r(p)$ denotes the closed disk of radius $r$ 
      centered at $p$.
      \item $A^o$ and $\partial A$ respectively denote the interior
        and boundary of $A$.
        \end{itemize}

As in Figure 5, define
  \begin{equation}
    O_{XY}=\bigg([-1,1]^2 -  \bigcup D^o_{2/3}(\pm 1,\pm1)\bigg)
    \times \{0\}.
\end{equation}
The union is taken over the $4$ disks we get by
taking all sign choices.
\begin{center}
  \resizebox{!}{2.7in}{\includegraphics{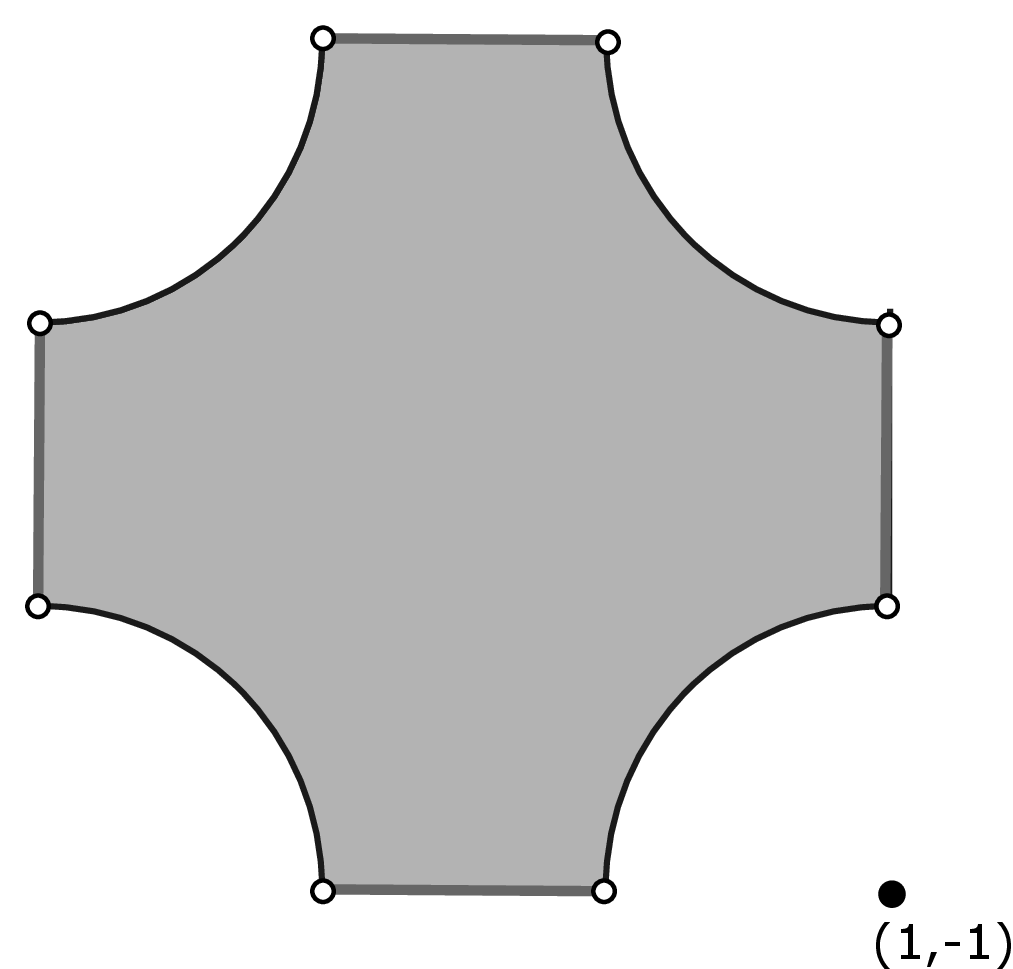}}
  \newline
  Figure 5:  The set $O_{XY}$
\end{center}
The set $O_{XY}$ is an octagon made from $8$ sides, with curved
and straight sides alternating.
The {\it octacross\/} is
  \begin{equation}
    O_{XYZ}=O_{XY} \cup O_{YZ} \cup O_{ZX}.
  \end{equation}
The octacross has $24$ edges, $12$ curved
and $12$ straight.  We call the curved edges {\it roads\/}.
The $12$ straight
edges cross each other in pairs, making $6$
{\it crosses\/}, two per coordinate axis.  If we
were to crush these crosses down to points,
the resulting graph would be the octahedral
graph.

\section{The Crossbridges}

Let
$$
A_{XY}=\bigg(D_{4/3}(1,-1) - D^o_{2/3}(1,-1)\bigg) \times \{0\},$$
\begin{equation}
B_{XY}=\partial D_1(1,-1) \times [-1/3,1/3].
\end{equation}
$A_{XY}$ and $B_{XY}$ are annuli which
intersect in the unit circle $\partial D_1(1,-1) \times \{0\}$.
One of the {\it crossbridges\/} is
\begin{equation}
  Q_{XY} = \{(x,y,t) \in A_{XY} \cup B_{XY}|\ \max(|x|,|y|) \geq 1\}.
\end{equation}
We are
chopping one quadrant out of $A_{XY} \cup B_{XY}$
to create $Q_{XY}$.
Figure 6 shows the projection of $Q_{XY}$ into the $XY$ plane.

\begin{center}
  \resizebox{!}{3.5in}{\includegraphics{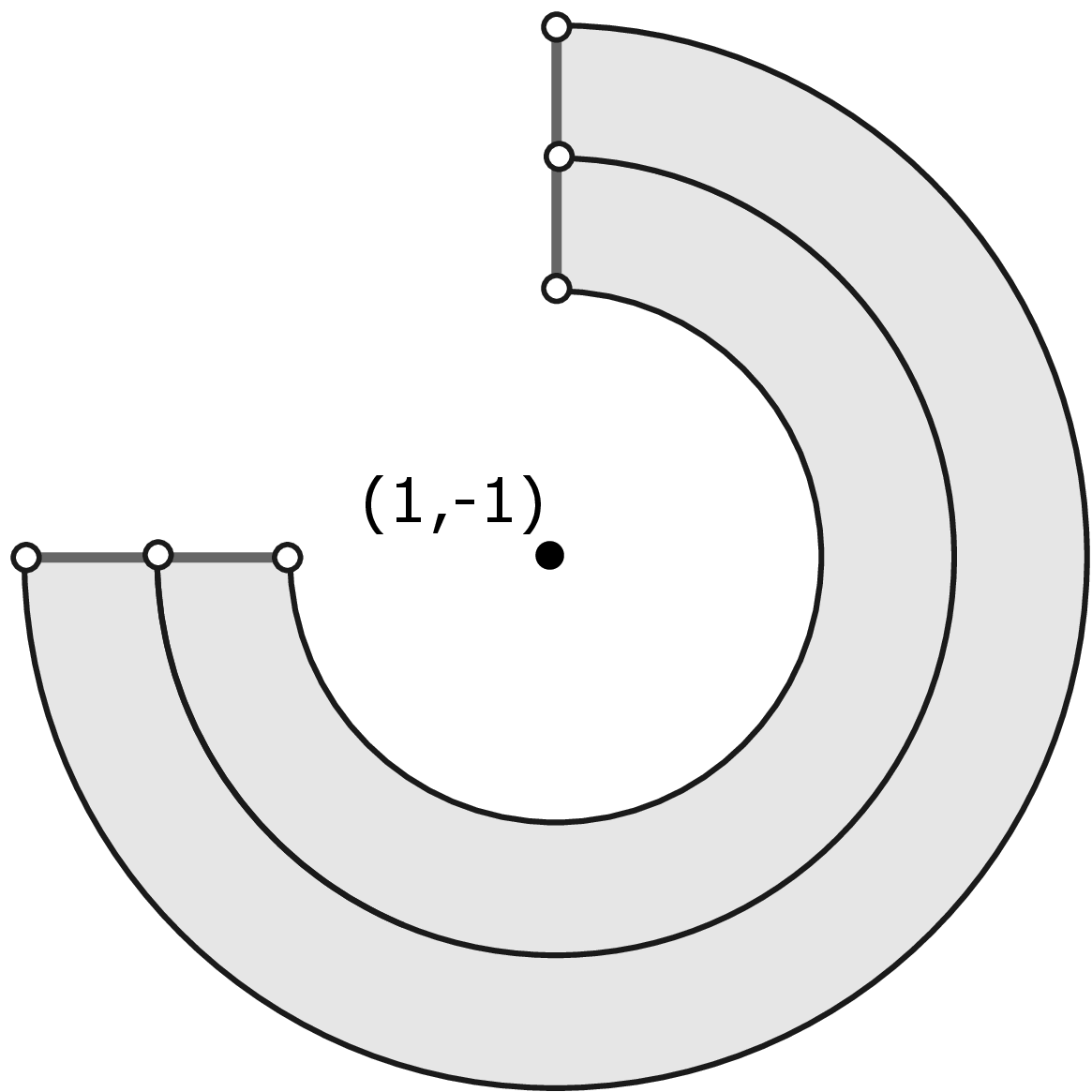}}
  \newline
  Figure 6:  The $XY$-projection of $Q_{XY}$.
  \end{center}

$Q_{XY}$ is made by sweeping out
a cross along three-quarters of a circle.
The boundary of $Q_{XY}$ consists of
$4$ curved edges  and $4$ straight
edges.  We call the curved edges {\it bridges\/}.
The straight edges in the boundary
cross in pairs, making $2$ crosses.
These crosses match $2$ of the crosses of the octacross.

The other crossbridges are
$Q_{YZ}$ and $Q_{ZX}$.

\section{Assembly}
\label{assemble}

Figure 7 shows
the $XY$-projection of the octacross $O_{XYZ}$ and the
crossbridge $Q_{XY}$.   This picture is a bit like
a head with one hooped earring.

\begin{center}
  \resizebox{!}{3.1in}{\includegraphics{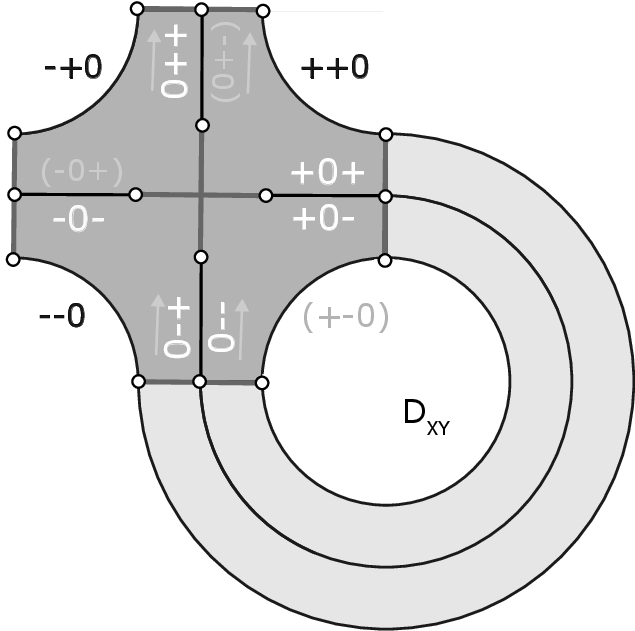}}
  \newline
  Figure 7:  The $XY$-projection of $O_{XYZ} \cup Q_{XY}$.
\end{center}

The set
\begin{equation}
  \label{head}
  H=O_{XYZ} \cup Q_{XY} \cup Q_{YZ} \cup Q_{ZX}
\end{equation}
is
like a head with $3$ hooped earrings.  Again, see Figure 0.
As Figure 7 suggests, $\partial H$ contains $3$
circles, each made from a road-bridge pair, one per coordinate plane.
We glue in the disks $D_{XY}$ and $D_{YZ}$ and $D_{ZX}$ bounded
by these circles.   That is, we define
\begin{equation}
  \label{Moebius}
  M=O_{XYZ} \cup Q_{XY} \cup Q_{YZ} \cup Q_{ZX} \cup  D_{XY} \cup D_{YZ}
  \cup D_{ZX}.
\end{equation}
$\partial M$ has $9=12-3$ roads and $9=12-3$ bridges.  Together
these make an embedded loop, a
combinatorial-length $18$ circuit in which the
roads and bridges alternate.  Combinatorially,   this
  is  the same circuit as the one in Figure 4.   For instance,
  looking at Figure 7 we can see that the roads
  $(+\!+\!0)$, $(-\!-\!0)$ and
  $(+0+)$, $(0\!-\!+)$ and
  $(+0-)$, $(0\!-\!-)$ connect across the shown crossbridge.
  This matches Figure 4. (Our vertical labels are meant
  to be read from bottom to top, as the arrows indicate.)

  Now I show that $M$ is an immersed Moebius band. 
        Let $M_{XY}$ denote the intersection of $M$ with the
        $XY$-plane.  Here $M_{XY}$ is a topological rectangle with
        alternating straight and curved edges.   See Figure 8.
  The closure of $M-M_{XY}-M_{YZ}-M_{ZX}$ is the union
  of the $B$-parts of the crossbridges,
  $B_{XY}$ and $B_{YZ}$ and $B_{ZX}$.
  Each of these pieces  is also a topological rectangle, with
  alternating straight and curved edges.
  We have a  cycle of $6$ pieces, with adjacent
  pieces sharing a common straight edge:
  \[
\begin{array}{ccccc}
M_{XY} & \leftrightarrow & B_{YZ} & \leftrightarrow & M_{ZX} \\[6pt]
\updownarrow &  &  &  & \updownarrow \\[6pt]
B_{ZX} & \leftrightarrow & M_{YZ} & \leftrightarrow & B_{XY}
\end{array}
\]

\begin{center}
  \resizebox{!}{3.7in}{\includegraphics{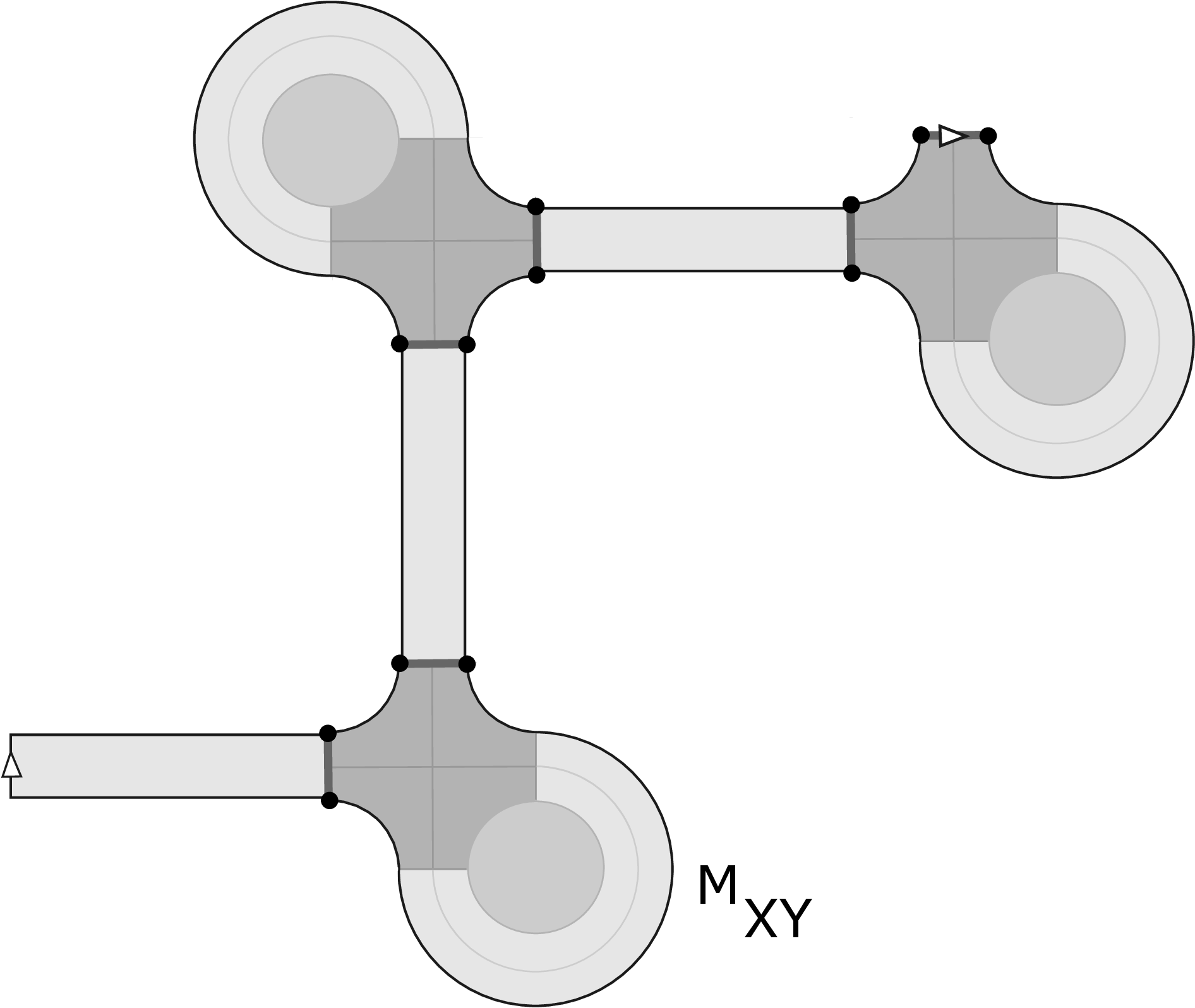}}
  \newline
  Figure 8:  $M_{XY}$ and the Moebius band $M^*$.
\end{center}

  We  treat each piece as a
  space in itself, and then form the identification space
  by gluing along the common straight boundary edges. This gives
  an abstract Moebius band $M^*$  that is decomposed into $6$
  topological rectangles. Figure 8 shows a fundamental domain for $M^*$.
The arrowed sides are meant to be glued.
    The obvious map $M^* \to M$ is an isometric immersion.

  \section{Adding the Pizza}
  \label{pizza}
  
  We work in the $3$-sphere $S^3=\R^3 \cup \infty$.
  Let $M$ be our immersed Moebius band, as in Equation \ref{Moebius}.
  As we mentioned in the introduction, $\partial M$ is cone-friendly,
meaning that each ray through the origin intersects $\partial M$ at
most once.
We  define the {\it pizza\/} as the cone of $\partial M$ to $\infty$:
    \begin{equation}
    P=\{tp|\ p \in \partial M,\ t \geq 1\} \cup \{\infty\}.
  \end{equation}
  We have
 $M \cap P = \partial M$.
 The union $M \cup P$ is Boy's surface.   We get Boy's surface
 by attaching an embedded disk to an immersed Moebius
 band along the topological circle $\partial M=M \cap P$.   Hence
 $M \cup P$ is an immersed projective plane.
 $M \cup P$ has the $6$-fold symmetry discussed at the end of
  \S \ref{octa}.

  We now look more closely at $P$, with a view towards making
  the cone-friendliness of $\partial M$ more transparent.
  $P$ is made from $18$ topological
  triangles, the {\it pizza slices\/},
  one per edge of $\partial M$.

\begin{center}
  \resizebox{!}{2.7in}{\includegraphics{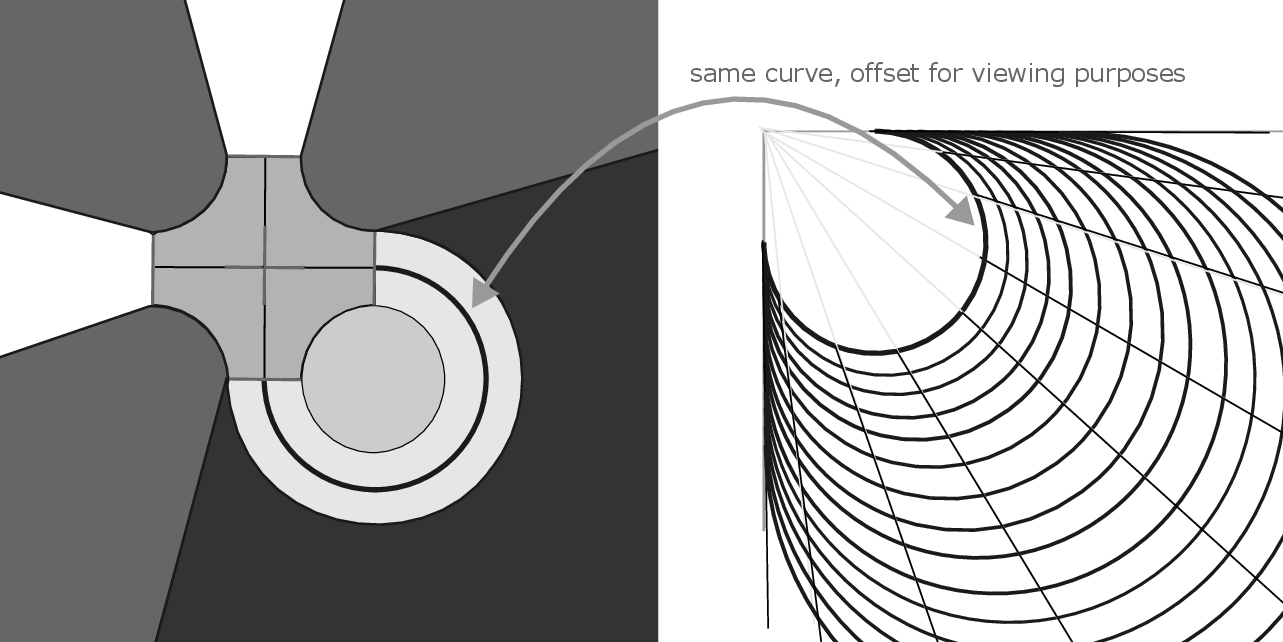}}
  \newline
  Figure 9:  The $4$ slices in the $XY$ plane and a hint of a
  non-planar slice
\end{center}

  Of these, $12$
lie in  coordinate
planes and $6$ do not. Each of the $6$ non-planar ones lies in an
elliptical cone, and each cone lies in a different orthant.
Figure 9 shows the $4$ pizza slices in the $XY$-plane.
Figure 9 also shows the non-planar slice contained in the $+-+$
orthant.  We are taking the double foliation of this slice by
(upward rising) rays and $(3/4)$-circles and projecting it down into the $XY$-plane.
The other pizza slices are obtained from these $5$ by the
symmetry discussed at the end of \S \ref{octa}.

\section{Discussion}
\label{airplane}

\noindent
{\bf The Airplane Movie:\/}
The airplane flies around and produces a set $H'$ which is the
same as $H$ (from Equation \ref{head})
topologically and close geometrically.
The movie then fills in the
three disks to make a set $M'$ which is likewise close to $M$.
Finally, the movie explicitly caps off $M'$
with embedded (curvy) disk $P' \subset \R^3$.
These notes differ from the movie in several ways.
\begin{itemize}
\item The airplane builds $H'$ in a different order than
  we build $H$.
  If you follow along the corridor in Figure 8, bottom
  to top, you are making the airplane's path.
\item The movie does not explain why $M'$ is a Moebius band
  or why $\partial M'$ is connected. (In general, nothing is explained
  or defined or discussed in the movie; everything is exhibited.)
  \item The disks $P'$ and $P$ differ in several ways.
  First, $P'$ lives in $\R^3$. This 
  breaks the $6$-fold symmetry discussed
  at the end of \S \ref{octa}.  Second, $P'$ seems to be smoother than
  $P$.  The cone point at $\infty$ for $P$ is rather complicated. 
\end{itemize}
I don't mean to suggest that these unexplained things are
deficiencies in the movie. The movie is amazing, and it has
different aims than these notes.
\newline
\newline
{\bf Intrinsic Flatness:\/}
The abstract Moebius band $M^*$ is a
  developable
  surface, as indicated by Figure 8.
  The universal cover $\widetilde M^*$ of $M^*$
  is an infinite
  topological strip.  The symmetries discussed at the
  end of \S \ref{octa} all correspond to symmetries of
  $\widetilde M^*$.  For instance,
  the order $3$ isometry of $\R^3$ which
  cycles the coordinates lifts to give a glide
  reflection of 
  $\widetilde M^*$.  The cube of this glide reflection is
  the generator of the deck transformation group.

  When we think of $\partial M^*$ as a curve in the
  boundary of an intrinsically flat Moebius band, it has
  total geodesic curvature $9\pi$.  In particular,
  $\partial M^*$ is not a paper Moebius band in the sense of
  [{\bf Sch\/}].  However, the situation here does raise an
  optimization question.   How much can you
  decrease the total geodesic curvature by applying an
  ambient homeomorphism of $\R^3$ which retains the
  intrinsic flatness of $M^*$ and also the cone-friendliness
  of $\partial M$?

  It is also worth mentioning that the intersection
  $P \cap \R^3$ is also
  intrinsically flat.  However, the geodesic curvature
  of $\partial P$, with respect to $P$, is rather complicated
  on the non-planar pizza slices. So, $M \cup P$ has a
  strange intrinsic geometry along parts of
$M \cap P$.

\section{Build Your Own}

   This section contains a kit 
   for building your own copy of the Moebius band $M$.
    The kit is set up in such a
    way that you don't need to glue in $3$ round disks.  These
    come already attached to the crossbridges.
    \newline
    \newline
    {\bf Step 0:\/} The next two pages have all the things you
    need to print out.  Print out 3 copies of each
    of the two pages.  I recommend using cardstock.
      \newline
    \newline
    {\bf Step 1A:\/} Cut out all $3$ copies of the
    piece which matches Figure 5.
    On the first copy, cut along all solid lines.
    On the second copy, cut only along the long solid line.
    Cut the third copy in half horizontally, 
    then cut on the vertical dotted lines.
    I found it useful to extend the cuts just a tiny bit,
    to make it easier to fit the pieces together all the way.
    \newline
    \newline
    {\bf Step 1B:\/} Take the two whole central pieces and
    slide them together along the cut slits.  Tape together.
    Now slot in the two halves of the central piece to
    finish the octacross.  Tape everything together as needed
    to make it  firm.
    \newline
    \newline
    {\bf Step 2A:\/} Cut out the $6$ strips and then cut along
    the dotted lines. Be sure to include the thin rectangles
    on the left side of each strip. I will explain these in the next
    step.
    Cut out the $3$ pieces made from a C-shaped annulus
    and the central disk, then cut along all dotted curves.
    (Note: Don't separate the central disk from the annulus.)
    The short crosscut on each annulus is just an
    access cut, to help you make the longer curved cut.
    \newline
    \newline
    {\bf Step 2B:\/} Slot two of the strips into one of the C shaped
    annuli along the cut slits, then tape together as needed.
    Finally, tape the access cut shut.   This makes one crossbridge
    together with the central disk already attached.  Now repeat
    for the other two crossbridges. One technical
    note:  The extra thin rectangle on the end of the strips is
    something like a tab.  The lengths of the strips are not precisely
    the same as the length of half the core curve, and this extra tab
    is supposed to absorb the mismatch.
    \newline
    \newline
    {\bf Step 3:\/} Join each crossbridge-plus-disk to the octacross
    as discussed in \S \ref{assemble}.  Tape together as needed.
    The disk part of the crossbridge-plus-disk should fit together
    with the relevant curved edge of the octacross.
    \newline

\clearpage
\thispagestyle{empty}
\noindent\includegraphics[
  width=\textwidth,height=\textheight,keepaspectratio,
]{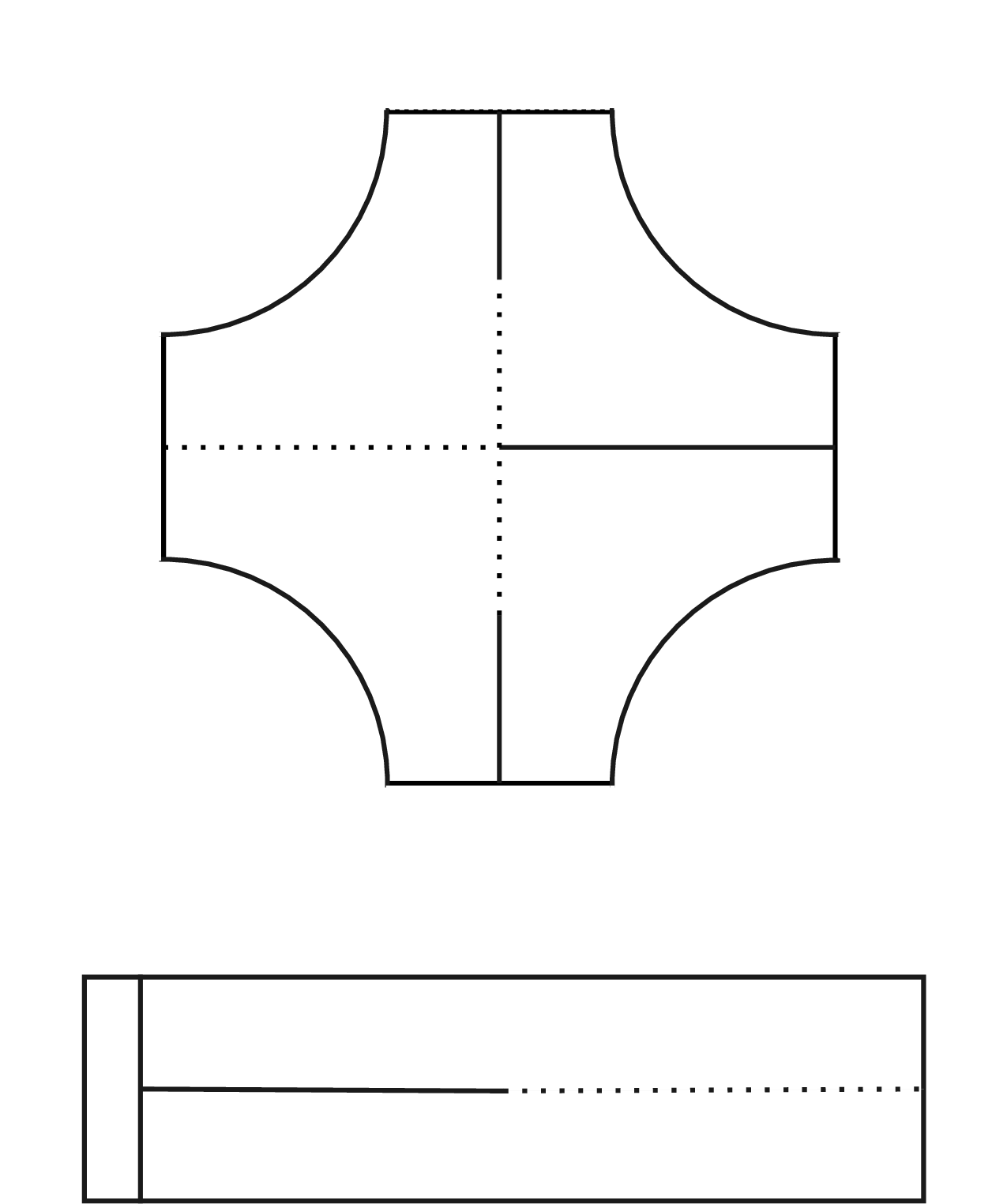}
\clearpage

\clearpage
\thispagestyle{empty}
\noindent\includegraphics[
  width=\textwidth,height=\textheight,keepaspectratio,
]{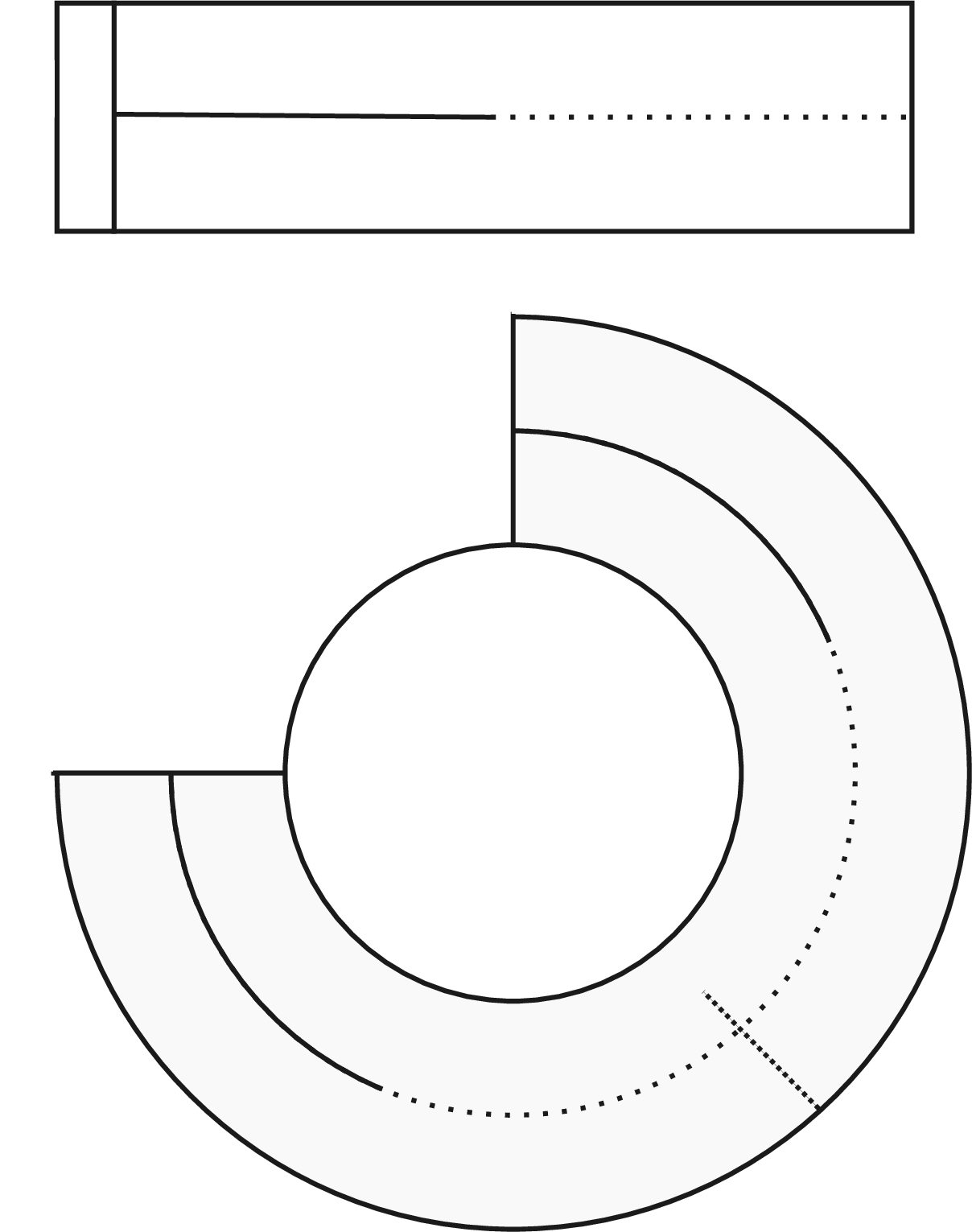}
\clearpage

  \section{A Rectilinear Model}
  \label{cube}

  Figure 10 shows a version of Figure 7 using alternate rectilinear
  models for the octacross, crossbridge, and disk.

\begin{center}
  \resizebox{!}{1.6in}{\includegraphics{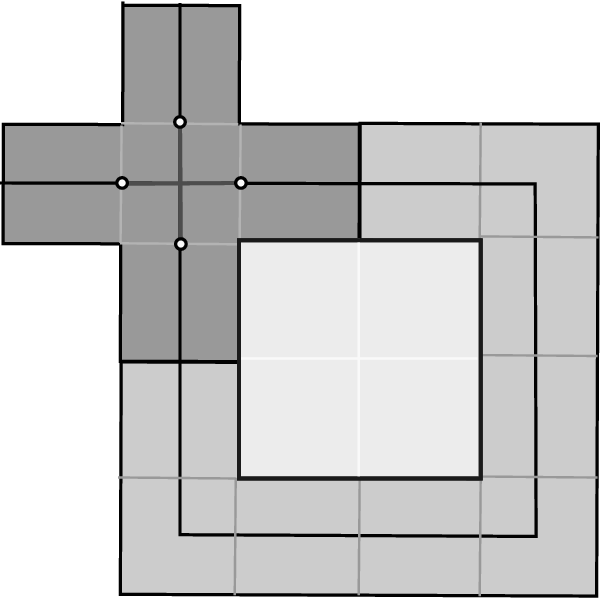}}
  \newline
  Figure 10:  Rectilinear models for the octacross, crossbridge, and disk.
\end{center}
The rectilinear pieces have especially nice thickenings:
\begin{itemize}
\item The octacross sits inside a union of $7$ cubes.
\item Each crossbridge sits inside a $C$-shaped chain of $9$ cubes.
\item Each square disk sits inside a union of $4=2 \times 2$ cubes.
\end{itemize}
The union $\Omega$ of these $7$ pieces is a topological ball
made from $46$ cubes.  Letting $M^{\square}$ denote our rectilinear
model for $M$, we have
$M^{\square} \subset \Omega$ and 
$\partial M^{\square} \subset \partial \Omega$.
We give an explicit description. Let
$Q_r$ denote the square of side-length $r$ centered at $(3,-3)$.
Scaling appropriately, we have
    \begin{equation}
\Omega=\Omega_{XY}\cup\Omega_{YZ}\cup\Omega_{ZX},
\qquad
M^{\square}=M^{\square}_{XY}\cup M^{\square}_{YZ}\cup M^{\square}_{ZX},
\end{equation}
\[
\Omega_{XY}=Q_4\times[-1,1],
\hskip 20 pt
M^{\square}_{XY}=(Q_4\times\{0\})\cup  (\partial Q_3 \times [-1,1]).
\]
    $\Omega$ is the union of three
    $8 \times 8 \times 2$ slabs, whose intersection
    is the $2 \times 2 \times 2$ cube centered at the origin.
    $\Omega$ looks like a cubical head with $3$ big ears.

    This analysis points out some flexibility in how we add
    the final disk to cap off Boy's surface.
    The complement $S^3-\Omega$ is also a ball, and
    $\partial M^{\square}$ is an embedded loop on the boundary.
    So, topologically speaking, it is easy to extend $\partial M^{\square}$
    to an embedded disk in $S^3-\Omega$.  The coning construction,
    which works in this model as well,
    does it in a canonical way.

    $\Omega$ is tiled by $368=46 \times 8$ unit cubes whose
    intersection
    with    $M^{\square}$  is a union of $1$, $2$, or $3$ adjacent
    faces.  Thus, you could build $(\Omega,M^{\square})$ by taking
    $368$ cubes, suitably painting their faces, and sticking them together.

\newpage

\section{References}
\noindent
[{\bf B\/}] W. Boy {\it \"Uber die Curvatura integra und die Topologie geschlossener Fl\"achen\/},
Math. Ann. Vol 57 (1903), 151--184.
\newline
\newline
[{\bf C\/}] M. Chas,  Instagram, \newline
www.instagram.com/reel/C-n6ExptzhI/?igsh=MTl0NGd4ZTg0YThveg
\newline
\newline
\noindent
[{\bf K\/}] R. Kirby {\it What is ... Boy's surface\/},
Notices of the A.M.S,  Vol 54, Number 10. (2007)
\newline
\newline
\noindent
[{\bf S\/}] Serbian Academy of Sciences,  Project
Ziva Mathematica, Youtube Video \newline
www.youtube.com/watch?v=9gRx66xKXek
\newline
\newline
[{\bf Sch\/}] R. E. Schwartz, {\it The Optimal Paper Moebius Band\/},
Annals of Mathematics, 2025.
\newline
\newline
\noindent
\end{document}